\documentclass{owrart}

\usepackage[all]{xy}

\usepackage{tikz}
\usetikzlibrary{arrows}
\usetikzlibrary{snakes}
\usetikzlibrary{shapes}
\newcommand{\Hom}{{\sf Hom }}\newcommand{\RHom}{\mathbf{R}{\sf Hom }}

\newcommand{\Ext}{{\sf Ext }}

\renewcommand{\mod}{{\sf mod \hspace{.02in}  }}
\newcommand{\fd}{{\sf fd \hspace{.02in}}}
\newcommand{\gr}{{\sf gr \hspace{.02in}}}
\newcommand{\per}{{\sf per \hspace{.02in}}}
\newcommand{\lten}{\overset{\mathbf{L}}{\otimes}}

\newtheorem{theorem}{Theorem}[section]
\newtheorem*{thm*}{Th{\'e}or{\`e}me}

\theoremstyle{remark}
\newtheorem{remark}[theorem]{Remark}
\newtheorem{example}[theorem]{Example}
\theoremstyle{definition}
\newtheorem{definition}[theorem]{Definition}

\begin{document}


\begin{talk}[O. Iyama, I. Reiten and S. Oppermann]{Claire Amiot}
{Preprojective algebras and Calabi-Yau duality}
{Amiot, Claire}

\noindent
The properties of the preprojective algebra are very different whether the associated quiver is of Dynkin type or not. However in both cases, one can construct from it a triangulated category of Calabi-Yau dimension $2$. In this note we explain the generalizations of this fact in the context of higher preprojective algebra, and we give some homological properties that characterize preprojective algebras.

\section{Classical case}

Let $k=\overline{k}$ be an algebraically closed field.
Let $Q$ be a finite quiver. The double quiver $\overline{Q}$ of $Q$ is defined from $Q$ by adding for each arrow $a\in Q_1$ an arrow $\overline{a}$ in the opposite direction. The \emph{preprojective algebra} of $Q$ is defined by \[\Pi_Q:=k\overline{Q}/\langle \prod_{a\in Q_1}[a,\overline{a}]\rangle.\] This notion has been defined by Gelfand and Ponomarev in \cite{GP}.

\begin{example}
Let $Q$ be the following quiver $\xymatrix{1\ar@(dr,ur)_x}$. Then we have $kQ\cong k[x]$. The preprojective algebra of $Q$ is presented by the quiver $\xymatrix{1\ar@(dr,ur)_x\ar@(dl,ul)^{\overline{x}}}$ with the relation $x\overline{x}-\overline{x}x=0$. That is $\Pi_Q\cong k[x,\overline{x}]$.
\end{example}

\begin{example}
Let $Q$ be the quiver $\xymatrix{1\ar[r]^a & 2}$. Then the preprojective algebra of $Q$ is presented by the quiver $\xymatrix{1\ar@/^/[r]^a & 2\ar@/^/[l]^{\overline{a}}}$ with the relations $a\overline{a}=\overline{a}a=0$.
\end{example}

To preprojective algebras we can associate some triangulated categories with a duality of Calabi-Yau type.

\begin{definition}[Kontsevich]
Let $\mathcal{T}$ be a $k$-linear triangulated category with finite dimensional Hom spaces. The category $\mathcal{T}$ is said to be $d$-\emph{Calabi-Yau} if for all objects $X,Y$ in $\mathcal{T}$, there exists an isomorphism $\Hom_\mathcal{T}(X,Y)\cong D\Hom_\mathcal{T}(Y,X[d])$ functorial in $X$ and $Y$.
\end{definition} 

The link between preprojective algebras and Calabi-Yau categories is given by the following classical result.

\begin{theorem}\cite{ES98,CB, BBK}\label{classical}
Let $Q$ be a finite quiver without oriented cycles. Then 
\begin{itemize}
\item if $Q$ is Dynkin, the preprojective algebra $\Pi_Q$ is finite dimensional selfinjective and the stable category $\underline{\mod}\Pi_Q$ is $2$-Calabi-Yau.
\item if $Q$ is not Dynkin, then $\Pi_Q$ is infinite dimensional and the bounded derived category $\mathcal{D}^{\rm b}(\fd \Pi_Q)$ of finite dimensional $\Pi_Q$-modules is $2$-Calabi-Yau.
\end{itemize}
\end{theorem}

The aim of this note is to explain how this result generalizes to higher preprojective algebras. 

\section{Higher preprojective algebras}

\begin{definition}[Iyama, Keller]
Let $\Lambda$ be a $k$-algebra of global dimension $d$. The ($(d+1)$-) preprojective algebra of $\Lambda$ is the tensor algebra $\Pi(\Lambda):={\sf T}_\Lambda \Ext^d_{\Lambda^{\sf e}}(\Lambda, \Lambda^{\rm e})$ where $\Lambda^{\rm e}=\Lambda^{\rm op}\otimes_k\Lambda$ is the envelopping algebra.
\end{definition}
As a tensor algebra, the preprojective algebra of $\Lambda$ is naturally equipped  with a  structure of positively graded algebra, with $\Lambda$ as degree $0$ subalgebra.

\begin{remark}
Ringel showed in \cite{Rin} (see also \cite{BGL}) that there is an isomorphism $\Pi_Q\cong \Pi(kQ)$.
\end{remark}

\begin{example} If $\Lambda$ is the the polynomial algebra in $d$ variables, it has global dimension $d$ and its associated preprojective algebra is isomorphic to the polynomial algebra in $d+1$ variables.
\end{example}

From now on $\Lambda$ we only consider the case where $\Lambda$ is finite dimensional. We denote by $\mathbb{S}_d^{-1}=-\lten_\Lambda \RHom_{\Lambda^{\rm e}}(\Lambda, \Lambda^{\rm e})[d]$ the autoequivalence of the derived category $\mathcal{D}^{\rm b}(\Lambda)$. 
We define the following subcategories of $\mathcal{D}^{\rm b}(\Lambda)$:
\[\mathcal{U}={\sf add}\{\mathbb{S}_d^{-p}\Lambda,\ p\in \mathbb{Z}\}\quad {\rm and}\quad \mathcal{U}^+={\sf add}\{\mathbb{S}_d^{-p}\Lambda,\ p\in \mathbb{N}\}.\]

\begin{definition}\cite{IO11,HIO12}
Let $\Lambda$ be a finite dimensional algebra of global dimension~$d$.  Then
\begin{itemize}
\item $\Lambda$ is said \emph{$d$-representation finite} (d-RF) if $\mathcal{U}=\mathcal{U}[d]$.
\item $\Lambda$ is said \emph{$d$-representation infinite} (d-RI) if $\mathcal{U}^+\subset \mod \Lambda$.
\end{itemize}
\end{definition}

Since the algebra $\Lambda$ has global dimension $d$, we have an isomorphism between ${\rm H}^0(\mathbb{S}_d^{-p}(\Lambda))$ and  $\Ext^d_{\Lambda^{\sf e}}(\Lambda, \Lambda^{\rm e})^{\otimes p}$ which is the homogeneous part of degree $p$ of $\Pi(\Lambda)$. Therefore one verifies that if $\Lambda$ is $d$-RF, then $\Pi(\Lambda)$ is finite dimensional, and if $\Lambda$ is $d$-RI, then $\Pi(\Lambda)$ is infinite dimensional.

For $d=1$, the functor $\mathbb{S}_1^{-1}$ is isomorphic to $\tau^{-1}$ the inverse of the Auslander-Reiten translation of the derived category $\mathcal{D}^{\rm b}(kQ)$. Then using Gabriel's theorem it is immediate to check the following equivalences: 
\[kQ \textrm{ is 1-RF } \Leftrightarrow Q \textrm{ is Dynkin }\Leftrightarrow kQ\textrm{ is representation-finite}\] \[kQ \textrm{ is 1-RI } \Leftrightarrow Q \textrm{ is non Dynkin }\Leftrightarrow kQ\textrm{ is representation-infinite}.\]

\begin{example}
Let $\Lambda$ be the algebra presented by the quiver $\xymatrix{1\ar[r]^a & 2\ar[r]^b & 3}$ with the relation $ba=0$. Then $\Lambda$ is $2$-RF. The preprojective algebra of $\Lambda$ is presented by the following quiver with relations: $$\xymatrix@R=.5cm@C=.5cm{&2\ar[dr]^b&\\1\ar[ur]^a & & 3\ar[ll]^c}\qquad ba=cb=ac=0.$$

There is a systematic way of constructing the preprojective algebra of an algebra of global dimension $2$. We refer to~\cite{Kel} for a description of the precise construction.
\end{example}

\begin{example}
Let $\Lambda$ be the algebra presented by the quiver $\xymatrix{1\ar[r]|y\ar@<1ex>[r]|x\ar@<-1ex>[r]|z & 2\ar[r]|y\ar@<1ex>[r]|x\ar@<-1ex>[r]|z & 3}$ with the commutativity relations. Then $\Lambda$ is $2$-RI. The preprojective algebra of $\Lambda$ is presented by the following quiver with the commutativity relations: $$\xymatrix@R=.7cm@C=.7cm{&2\ar[dr]\ar@<-.5ex>[dr]\ar@<.5ex>[dr]&\\1\ar[ur]\ar@<-.5ex>[ur]\ar@<.5ex>[ur] & & 3\ar[ll]\ar@<.5ex>[ll]\ar@<-.5ex>[ll]}.$$ 
\end{example}

\section{Results}

\subsection{$d$-RI case}
The first result can be understood as a generalization of the non Dynkin case of Theorem \ref{classical}. Furthermore it gives a homological characterization for the preprojective algebras of $d$-RI algebras. 

\begin{theorem}\cite{Kel, HIO12, MM, AIR11}\label{dRI}
Let $\Gamma=\bigoplus_{i\geq 0}\Gamma_i$ be a graded algebra with finite dimensional degree zero part $\Lambda:=\Gamma_0$.
Then the following are equivalent 
\begin{enumerate}
\item $\Lambda$ is $d$-RI, has global dimension $d$ and $\Gamma\cong \Pi(\Lambda)$ as graded algebras;
\item $\Gamma$ is $(1)$-twisted bimodule $(d+1)$-Calabi-Yau.
\end{enumerate}
\end{theorem}

Property (2) is an algebraic (and graded) enhancement of the $(d+1)$-Calabi-Yau property for the category $\mathcal{D}^{\rm b}(\fd\Gamma)$. This is the bimodule property that must be satisfied by $\Gamma$ to ensure that $\mathcal{D}^{\rm b}(\fd\Gamma)$ is $(d+1)$-Calabi-Yau. It is given as follows: 
\[(2)\Leftrightarrow \Gamma\in \per \Gamma^{\rm e}\ \textrm{and}\quad \RHom_{\Gamma^{\rm e}}(\Gamma, \Gamma^{\rm e})[d+1]\cong \Gamma(1)\ \textrm{in }\mathcal{D}(\gr \Gamma^{\rm e}).\]
Here $\Gamma(1)$ is the graded bimodule whose degree $p$ part is $\Gamma_{p+1}$. 

Theorem \ref{dRI} is stated in \cite[Thm 4.35]{HIO12}.  The implication $(1)\Rightarrow (2)$ follows from \cite[Thm 4.8]{Kel}, while $(2)\Rightarrow (1)$ was shown indenpendently in \cite[Thm 4.8]{MM} and \cite[Thm 3.4]{AIR11}.

\subsection{$d$-RF case}

The next result is the $d$-RF analogue of Theorem \ref{dRI}. It can be seen as a generalization of the Dynkin case of Theorem \ref{classical}, and gives a homological characterization of the preprojective algebras of $d$-RF algebras. 
\begin{theorem}\label{dRF}\cite{Dugas, AO14}
Let $\Gamma=\bigoplus_{i\geq 0}\Gamma_i$ be a finite dimensional graded algebra. Denote by $\Lambda$ its degree zero part.
Then the following are equivalent 
\begin{enumerate}
\item $\Lambda$ is $d$-RF, has global dimension $d$ and $\Gamma\cong \Pi(\Lambda)$ as graded algebras;
\item $\Gamma$ is selfinjective and $(1)$-twisted stably bimodule $(d+1)$-Calabi-Yau.
\end{enumerate}
\end{theorem}
 Here again, property (2) is an algebraic (and graded) enhancement of the $(d+1)$-Calabi-Yau property for the category $\underline{\mod}\Gamma$. This is the bimodule property that must be satisfied by $\Gamma$ to ensure that $\underline{\mod}\Gamma$ is $(d+1)$-Calabi-Yau. It is given as follows: 
\[(2)\Leftrightarrow  \Hom_{\Gamma^{\rm e}}(\Gamma, \Gamma^{\rm e})[d+2]\cong \Gamma(1)\ \textrm{in }\underline{\gr} \Gamma^{\rm e}.\]

The implication $(1)\Rightarrow (2)$ is shown in \cite[Thm 3.2]{Dugas}, while the implication $(2)\Rightarrow (1)$ is shown in \cite[Thm 3.1]{AO14}.

\subsection{Beyond the RF/RI cases}

Most of algebras of global dimension $d\geq 2$ are neither $d$-RF nor $d$-RI. So one could ask how Theorems \ref{dRI} and \ref{dRF} can be extended to general preprojective algebras. Here we focus on the case where the preprojective algebra is finite dimensional. 

In general the finite dimensional preprojective algebras are not selfinjective but their behaviour is still similar to the one of the preprojective algebras of $d$-RF algebras. In the case $d=2$, it is shown in \cite{A08} that when $\Pi(\Lambda)$ is finite dimensional, it is the endomorphism ring of a cluster-tilting object in a certain $2$-Calabi-Yau category. Keller and Reiten proved that such algebras are Gorenstein (that is ${\rm proj dim}D\Pi={\rm inj dim }\Pi<\infty$). Hence the correct analogue Calabi-Yau triangulated category is given by the stable category of maximal Cohen-Macaulay $\Pi$-modules. Indeed they proved in \cite{KR07} that the category $\underline{\sf CM}\;\Pi(\Lambda)$ is $3$-Calabi-Yau.

These results were the motivation for the following characterization of finite dimensional preprojective algebras.
\begin{theorem}\cite{AO14}\label{CM}
Let $\Gamma=\bigoplus_{i\geq 0}\Gamma_i$ be a (non trivially) graded finite dimensional algebra. Denote by $\Lambda$ its degree zero part. Assume that 
\begin{itemize}
\item[(a)]$\Gamma$ is Gorenstein of dimension $\leq d-1$;
\item[(b)] there is an isomorphism $\RHom_{\Gamma^{\rm e}}(\Gamma,\Gamma^{\rm e})[d+2]\cong \Gamma(1)$ in $\mathcal{D}^{\rm b}(\gr \Gamma^{\rm e})/\per\gr \Gamma^{\rm e}$.
\item[(c)] $\Ext^{i}_{\Gamma^{\rm e}}(\Gamma,\Gamma^{\rm e}(j))=0$ for any $i\geq 1$ and any $j\leq -1$.
\end{itemize}
Then $\Lambda$ has global dimension $d$ and $\Gamma\cong \Pi(\Lambda)$ as graded algebras.
\end{theorem}

Here property $(b)$ is again an algebraic (and graded) enhancement of the $(d+1)$-Calabi-Yau property of the category  $\underline{\sf CM}\;\Gamma$. 

One also shows in \cite{AO14} that these properties are satsified by finite dimensional preprojective algebras in the case $d=2$ and $d=3$ using the description of $\Pi(\Lambda)$ in term of quivers with relations.

\end{talk}

\end{document}